\documentclass{ifacconf}
\usepackage{graphicx}      
\usepackage{natbib}        
\usepackage{graphicx}
\usepackage{amsmath}
\usepackage{tabularx} 
\usepackage{graphicx} 
\usepackage{float}
\usepackage{multirow}
\usepackage{amssymb}

\newcommand{\cC}{\mathcal{C}}

\newcommand{\R}{\mathbb{R}}

\newcommand{\cL}{\mathcal{L}}

\newcommand{\cJ}{\mathcal{J}}

\usepackage{color}

\newcommand{\bU}{\mathbf{U}}
\newcommand{\bs}{\mathbf{s}}

\newcommand{\bW}{\mathbf{W}}
\newcommand{\bb}{\mathbf{b}}

\newcommand{\bo}{\mathbf{0}}

\newcommand{\bA}{\mathbf{A}}
\newcommand{\bQ}{\mathbf{Q}}
\newcommand{\bB}{\mathbf{B}}
\newcommand{\bR}{\mathbf{R}}

\newcommand{\bx}{\mathbf{x}}
\newcommand{\by}{\mathbf{y}}
\newcommand{\bu}{\mathbf{u}}
\newcommand{\bff}{\mathbf{f}}

\begin{document}
\begin{frontmatter}

\title{Supervised learning for kinetic consensus control} 


\author[First]{Giacomo Albi} 
\author[Second]{Sara Bicego} 
\author[Second]{Dante Kalise}

\date{\small{\textit{${}^*$ Department of Computer Science, University of Verona\\ Strada le Grazie 15 - 37134 Verona, Italy (email: giacomo.albi@univr.it)\\
${}^{**}$ Department of Mathematics, Imperial College London, South Kensington \\ Campus - SW72AZ London, UK (email: s.bicego21@imperial.ac.uk)\\ 
${}^{***}$ Department of Mathematics, Imperial College London, South Kensington \\ Campus - SW72AZ London, UK (email: d.kalise-balza@imperial.ac.uk)
}}}

\address[First]{Department of Computer Science, University of Verona, Strada le Grazie 15 - 37134 Verona, Italy (email: giacomo.albi@univr.it).}
\address[Second]{Department of Mathematics, Imperial College London, South Kensington Campus - SW72AZ London, UK (email: {s.bicego21,dkaliseb}@imperial.ac.uk)}

\begin{abstract}                
In this paper, how to successfully and efficiently condition a target population of agents towards consensus is discussed. To overcome the curse of dimensionality, the mean field formulation of the consensus control problem is considered. Although such formulation is designed to be independent of the number of agents, it is feasible to solve only for moderate intrinsic dimensions of the agents space. For this reason, the solution is approached by means of a Boltzmann procedure, i.e. quasi-invariant limit of controlled binary interactions as approximation of the mean field PDE. The need for an efficient solver for the binary interaction control problem  motivates the use of a supervised learning approach to encode a binary feedback map to be sampled at a very high rate. A gradient augmented feedforward neural network for the Value function of the binary control problem is considered and compared with direct approximation of the feedback law.
\end{abstract}

\begin{keyword}
Multi-agent systems, optimal feedback control, mean field models, supervised learning, opinion dynamics.
\end{keyword}

\end{frontmatter}

\section{Introduction}

Social behaviours can be seen as the result of a suitable combination of endogenous population interactions and external influences. How to successfully condition a population of agents towards a designed purpose is a fascinating question, whose answer is being widely researched \citep{li_survey_2019}. 

The formulation of such a problem in a dynamic optimization framework ensures the
availability of control synthesis methods, which nonetheless come with the huge drawback
of the curse of dimensionality. The problem reads as the minimization of a cost
functional subject to individual-based interaction dynamics, thus its solution
easily becomes unfeasible to compute as the number of agents in the population
grows. The natural way of circumventing this is using a multiscale approach working with the population density instead of its microscopic state. For a number of $N\to\infty$ interacting
agents, this leads to a mean field formulation of the control problem.
Although mean field optimal control problems are designed to be independent
of the number of agents, they are computationally feasible only for moderate
intrinsic dimensions $d$ of the agents’ state space. For this reason, we rely on the approximation of the mean field PDE governing the evolution of probability distribution characterizing the agents' population. An approximate solution is obtained as a result of a Boltzmann type dynamics \citep{albi_boltzmann_2017,albi_binary_2013,albi_mean_2017}.  This procedure provides the approximated mean field solution as limit of a reduced problem, modeling the interactions taking place only within controlled couples of agents. { We refer the reader to \citep{suboptimal_mf_2018} for a generalization of this procedure when allowing only a subset of agents in the system to be influenced by an external control signal, and to \citep{local_action} for a system-generalized control action.}

The efficiency of this Boltzmann approach is linked to the availability of a sufficiently fast solver for a binary interaction control problem, that is, an optimal control problem for a reduced system of two agents, which is sampled at a very high-frequency rate. We address this computational requirement by enconding a binary feedback map  by means of a supervised learning procedure \citep{kang_algorithms_2020,darbon_overcoming_2020,azmi_optimal_2021,albi_gradient-augmented_2022},  which is trained upon synthetic data from sampling a feedback law generated by a state-dependent Riccati equation approach (SDRE) \citep{C97,BLT07,Astolfi2020}. This feedback law corresponds to an approximation of the associated optimal feedback law characterized by the solution of a Hamilton-Jacobi-Bellman PDE. Despite being suboptimal, the SDRE law locally asymptotically stabilizes the dynamics, and can be easily computed by the sequential solution of algebraic Riccati equations, providing a reasonable alternative in high-dimensional settings where the numerical approximation of optimal feedback laws is prohibitively expensive.

 The rest of the paper is organized as follows. In Section 2 we introduce the mean field formulation of the addressed consensus problem, and in Section 3 we present a consistent alternative description of Boltzmann type. In Section 4 the state-dependent Riccati equation approach is presented, and its numerical approximation through supervised learning is discussed in  Section 5. A computational assessment dealing with control of first order opinion dynamics can be found in Section 6.

\section{Mean Field Consensus Problem}
We consider a population of $N_a$ agents evolving according to interaction dynamics of form:
\begin{equation}\label{CStest_dynamics}
   \dot{x}_i = \frac{1}{N_a} \sum\limits_{j = 1}^{N_a} P(x_i,x_j)(x_j - x_i) + u_i\,\quad x_i(0)=x_i^0\,,
\end{equation}
where the kernel $P(x_i,x_j)$ models the communication between agents with states $x_i\in\R^d$, and the control variable $\bu = (u_1,...,u_{N_a})$ aims at steering the system towards a consensus state $\Bar{x}=\dfrac1{N_{a}}\sum\limits_{i=1}^{N_a} x_i$. We express this goal as an infinite horizon nonlinear stabilization problem 
\begin{equation}\label{CStest_cost}
    \underset{\bu(\cdot) \in \mathcal{L}^2(\R_+ ; \mathbb{R}^{d \times N_a})}{\min}  \int\limits_{0}^{\infty} \frac{1}{N_a} \sum_{i=1}^{N_a} \|x_i - \Bar{x}\|^2 + \beta||u_i||_2^2 dt\,,
\end{equation}
subject to \eqref{CStest_dynamics}. A natural feature of agent-based models is that the number of interacting agents can become prohibitively large. Hence, as the number of agents $N_a$ grows, instead studying the microscopic, individual-based optimal control problem \eqref{CStest_dynamics}-\eqref{CStest_cost}, one can conveniently model the population by means of the density distribution of agents
\begin{equation}
    f = f(t;x), \qquad t\geq0,\qquad x \in \R^d\,,
\end{equation}
which evolves in time according with dynamics of the form
\begin{equation}\label{mf}
    \partial_t f = - \nabla_x\cdot \bigg[( \mathcal{P}[f] + u)f\bigg]\,,
\end{equation}
where the mean field interaction force $\mathcal{P}$ relative to the distribution $f$ reads
\begin{equation}
    \mathcal{P}[f(x)] = \int\limits_{\Omega} P(x,s)(s-x) f(s) ds\,.
\end{equation}
The optimal solution of the mean field optimal control problem -- obtained as combination of \eqref{mf} with a suitable cost functional -- is, by construction, independent of the number of agents, since it models the macroscopic behaviour of the population as a whole. However, the mean field optimal control solutions are meant to be computed via first-order optimality conditions, whose complexity is linked to the dimensionality $d$ of the state space: even for moderate values of $d$, the computational cost can be formidably high \citep{bensoussan_mean_2013}\citep{fornasier_mean-field_2014}.

\section{Boltzmann-type Formulation}
{ To circumvent the difficulties related to the solution of the mean field control problem, here
we aim at modeling the evolution in time of the population density function $f(t,x)$ from a kinetic viewpoint instead}. To this end,  we assume two agents with states $x_i,x_j\in\R^d$ interacting  according to the binary rule
\begin{equation}\label{post_int_dyn}
	\begin{aligned}
		x_i^* &= x_i + \eta\,\bigg( P(x_i,x_j)(x_j-x_i) + \, u(x_i,x_j)\bigg)\\
		x_j^* &= x_j + \eta\,\bigg( P(x_j,x_i)(x_i-x_j) +\, u(x_j,x_i)\bigg)\,,
	\end{aligned}
\end{equation}
where $\eta$ measures the strength of the interaction,  and $(x_i^*,x_j^*)$ are the post-interaction states. Hence, 
the evolution of $f(t,x)$ is driven by a Boltzmann-type dynamics:
\begin{equation}\label{boltzmann}
	\partial_t f(t,x)  = \lambda  \mathcal{Q}_{\eta,u}(f,f) (t,x)\,,
\end{equation}
where $\lambda $ is a parameter describing the  interaction frequency, and the operator $Q_{\eta,u}(f,f)$ accounts for the gain and loss of particles located a certain position $x$ at time $t$, as follows
	\begin{align}
		\mathcal{Q}_{\eta,u}(f,f) = \mathcal{Q}^{+}_{\eta,u}(f,f) - \mathcal{Q}^{-}_{\eta,u}(f,f)  
	\end{align}
with 
	\begin{align*}
		\mathcal{Q}_{\eta,u}^+(f,f)(t,x) &= \int\limits_{\Omega}\dfrac{1}{\mathcal{J}_\eta} f(t,{}^*x)f(t,{}^*s) d\bs\,,\\
		\mathcal{Q}_{\eta,u}^-(f,f)(t,x) &=f(t,x) \int\limits_{\Omega}f(t,s)ds\,,
	\end{align*}
and where $({}^*x_i,{}^*x_j)\longmapsto(x_i,x_j)$ are the pre-interaction states associated to \eqref{post_int_dyn},  and $\mathcal{J}_\eta$ represents the Jacobian of the binary interaction \eqref{post_int_dyn}. { The interest in solving \eqref{boltzmann}, arises when considering under a quasi-invariant scaling (i.e. $\eta = \varepsilon,\lambda=\varepsilon^{-1}$), as this provides us with the following consistency theorem between the mean field evolution of the dynamics and their Boltzmann formulation. We refer the reader to \citep{albi_mean_2017} for detailed derivation and proof of the result.
\paragraph*{Theorem 1}
Let $\eta\geq0$, $\varepsilon>0$, $P(\cdot,\cdot)\in\mathcal{L}^2_{loc}$ at all times $t\in[0,+\infty)$, and we consider a weak solution $f$ of \eqref{boltzmann} from initial condition $f_0(x)$. Furthermore, we introduce the scaling $\eta = \varepsilon,\,\lambda=\varepsilon^{-1}$ for the binary interaction rule, and we define $f^\varepsilon(t;x)$ to be a solution for the associated scaled version of \eqref{boltzmann}. Then, as $\varepsilon\to0$, we have pointwise convergence (up to subsequences) of the scaled solution $f^\varepsilon(t;x)$ to the solution $f(t;x)$ of \eqref{mf}.}

Different numerical schemes  can be derived to simulate the kinetic dynamics, \citep{albi_binary_2013}. In particular, the evolution of $f=f(t,x)$ can be approximated by means of Direct Simulation Monte Carlo Methods, introducing a forward Euler discretization as follows
\begin{align}\label{euler}
	f_{n+1} &= f_{n} + \Delta t \lambda\bigg(\mathcal{Q}^+_{\eta,u}(f_n,f_n) - \mathcal{Q}^-_{\eta,u}(f_n,f_n)\bigg)\\
	&= \big(1 - \Delta t\lambda\big)f_{n} + \Delta t\lambda\, Q_{\eta,u}^{+}(f_n,f_n),
\end{align}
with $\Delta t\leq \varepsilon$ to preserve positivity of the solution $f_{n+1}$. Thus, sampling $N_{\rm sample}$ particles from the initial distribution $f_0(x)= f(0,x)$ we can approximate the solution of \eqref{euler} via stochastic simulation of the binary interaction \eqref{post_int_dyn}. 

The convenience of this Boltzmann-type description relies on the possibility of approximating the behaviour of the population as the quasi-invariant limit of binary interactions, meaning that at each time step the agents are influenced only within couples. This heavily tackles down the computational complexity involved, since we are now considering the combination of many $2-$agents subproblems. The number of interacting couples depends on the frequency parameter $\lambda = 1/\varepsilon$: a choice $\varepsilon\ll1$ leads to weak, but frequent interactions, which is the typical case of mean-field models. Nonetheless, this requires an efficient solver for the reduced consensus problem.

\newpage
\section{State Dependent Riccati Equation}
{In this section, we aim at solving the reduced 2-agents problem, for which the states --encoding the position of both the coupled agents $i$-$j$ -- are denoted as a single variable $\bx(t)=(x_i(t),x_j(t))^{\top}\in\R^{2d}$. Similarly, we use bold notation when referring to the interaction force and the control variable associated to the dynamical system for $\bx$.}

The binary consensus problem resulting from the microscopic formulation\eqref{CStest_dynamics}-\eqref{CStest_cost} can be written as a nonlinear quadratic regulator problem (NLQR) 
\begin{equation}\label{ocp}
	\underset{\bu(\cdot)\in \bU}{\min}\cJ_{\bx_0}(\bu(\cdot)):=\int\limits_0^\infty \bx^{\top}\!(s)\bQ\bx(s)\,+\,\bu^{\top}\!(s)\bR\bu(s)\,ds\,,
\end{equation}
subject to nonlinear, control-affine dynamics
\begin{align}\label{dynamics}
	\dot \bx(t)=\bff(\bx(t))+\bB\bu(t)\,,\quad\bx(0)=\bx_0\,,
\end{align}
where $\bu(\cdot)\in\bU=\{\bu(t):\, \R^+\rightarrow \R^{2d}, \text{measurable}\}$ is an unbounded control variable, $\bQ\in\R^{2d\times 2d}$ is a symmetric positive semidefinite matrix, and  $\bR\in\R^{2d\times 2d}$ is symmetric positive definite.  The control operator $\bB:\R^{2d\times 2d}$, and  the system dynamics $\bff(\bx):\R^{2d}\rightarrow\R^{2d}$ are $\cC^1(\R^{2d})$ and such that $\bff(\bo)=\bo$ and $\bB(\bo)=\bo$.
Using dynamic programming, the optimal feedback law $\bu(\cdot)$ solving \eqref{ocp} can be recovered in terms of the value function of the control problem
\begin{equation}
	V(\bx)=\underset{\bu(\cdot)\in \bU}{\inf}\cJ_{\bx}(\bu(\cdot))\,,
\end{equation}
solving the following first-order, static, nonlinear Hamilton-Jacobi-Bellman PDE
\begin{align}\label{hjb}
	\nabla V(\bx)^{\top}\!\bff(\bx)&-\frac14\nabla V(\bx)^{\top}\!\bW(\bx)\nabla V(\bx)\!+\!\bx^{\top}\!\bQ\bx=0\,,
\end{align}
where $\bW=\bB\bR^{-1}\bB^{\top}$. Once the function $V(\bx)$ is computed, the associated optimal feedback is given by
\begin{equation}\label{feedback}
	\bu(\bx)=-\frac{1}{2}\bR^{-1} \bB^{\top} \nabla V(\bx)\,.
\end{equation}
Solving \eqref{hjb} can be in general difficult and expensive from a computational point of view. The value function $V(\cdot)$ maps variables living in $\R^{2d}$, where the dimension $d$ can be arbitrarily high. Equation \eqref{hjb} is a nonlinear PDE, thus it can be unfeasible to solve via standard methods even for moderate dimensional optimal control problems ($d>3$).

\subsection{Algebraic Riccati Equation and state-dependence}
In a simplified setting, under further assumptions of linearity for the free dynamics $\bff(\bx)=\bA\bx$ and making the ansatz $V(\bx)=\bx^{\top}\!\Pi\bx$ with $\Pi\in\R^{2d\times 2d}$, the optimality condition \eqref{feedback} can be written as 
\begin{equation}\label{feedbackare}
	\bu(\bx)=-\bR^{-1} \bB^{\top} \Pi\bx\,,
\end{equation}
where $\Pi$ is a positive definite solution of the Algebraic Riccati Equation (ARE)
\begin{equation}\label{are}
	\bA^{\top}\Pi+\Pi \bA-\Pi \bB\bR^{-1}\bB^{\top}\Pi+\bQ=0\,.
\end{equation}

Even though the class of problems being addressed in this paper is characterized by non-linearity in the free dynamics, a similar approach arises when casting \eqref{dynamics} in semilinear form:
\begin{equation}\label{semilinear}
    \bff(\bx) = \bA(\bx)\bx\,,\qquad \dot{\bx} =  \bA(\bx)\bx + \bB\bu\,. 
\end{equation}
In this setting, the solution of the HJB PDE \eqref{hjb} can be approximated by a state-dependent Riccati equation (SDRE)
\begin{align}\label{sdre}
	\bA^{\top}\!(\bx)\Pi(\bx)&\!+\!\Pi(\bx) \bA(\bx)\!-\!\Pi(\bx) \bW\Pi(\bx)\!+\!\bQ\!=\!0\,.
\end{align}

In the linear case, the ARE \eqref{are} directly comes from the HJB PDE \eqref{hjb} by considering the ansatz $V(\bx)=\bx^{\top}\Pi\bx$ for the associated value function. Thus, the feedback \eqref{feedbackare} resulting from the ARE solution $\Pi$ coincides with the optimal control variable resulting from \eqref{feedback}. The consistency between the ARE and the Dynamic Programming solutions is not readily available when dealing with nonlinear dynamics of the form \eqref{semilinear}. This is due to the state-dependence in the SDRE solution $\Pi(\bx)$, which leads to 
\begin{equation}\label{ansatz}
\begin{aligned}
    V(\bx)&=\bx^{\top}\Pi(\bx)\bx\,,\\ \nabla V(\bx)&=2\Pi(\bx)\bx + \varphi(\bx)\,,
\end{aligned}
\end{equation}
where $\varphi(\bx)$ is a $2d$-dimensional vector-valued function such that
\begin{equation}
    \varphi_k(\bx) = \sum\limits_{i,j = 1}^{2d} x_i x_j \dfrac{\partial \Pi(\bx)_{i,j}}{\partial x_k}\,.
\end{equation}
Thus, when substituting \eqref{ansatz} in the HJB PDE \eqref{hjb}, we do not recover the SDRE \eqref{sdre}, due to the presence of an additional term associated to the $\varphi(\bx)$ component in $\nabla V(\bx)$. For this reason, the feedback law 
\begin{equation}\label{feedbacksdre}
    \bu(\bx)=-\bR^{-1} \bB^{\top} \Pi(\bx)\bx\,,
\end{equation}
is a suboptimal approximation to the HJB feedback. Nevertheless, under stabilizability assumptions, the SDRE feedback law is locally asymptotically stabilizing \citep{BLT07}.

\subsection{Freezing coefficients in the Riccati Equation}
The main computational bottleneck associated to the synthesis of the SDRE feedback law is that eq. \eqref{sdre} cannot be solved analytically for  for $\Pi(\bx)$, and needs to be realized in a model predictive control fashion along a trajectory, as proposed in \citep{BLT07}. Given the current state $\bx$ of the system, we assume the operator $\Pi(\bx)$ to be a positive definite matrix in $\Pi\in\R^{d\times d}$, meaning that \eqref{sdre} reduces to its algebraic form \eqref{are}.  

This procedure can be useful to generate suboptimal approximations of the controlled trajectories associated to infinite horizon control problems of the form \eqref{ocp}-\eqref{dynamics}. While evolving along a trajectory, we assume the system to be in a configuration $\bar\bx$. By freezing every state-dependent operator accordingly with the current state $\bar\bx$, we obtain an ARE to be solved for the frozen SDRE operator $\Pi(\bar\bx)$, and the associated feedback law $\bu(\bar\bx)$ can be recovered via \eqref{feedbacksdre}. Then, we let the system evolve according with the $\bu(\bar\bx)$-controlled dynamics for a short time frame, after which the procedure is repeated. 

 Even if we assume that this SDRE approach generates asymptotically stable closed-loop solutions \citep{BLT07}, a main limitation persists, residing in the implementation of a sufficiently efficient ARE solver to enable a high-frequency sampling of controlled binary interactions \eqref{post_int_dyn}. For this efficiency purpose, we rely on supervised learning approximation models \citep{kang_algorithms_2020}, \citep{darbon_overcoming_2020} to encode the control action in a  neural network.

\section{Supervised Learning Approximation}
We populate a training set for the control law by solving in an offline phase the frozen SDREs  for a collection of states associated to $N_s$ sampled couples in $\R^d\times\R^d$ of interacting agents $\mathcal{X}_t=\{(\bx=x_i,x_j)^{(k)}\}_{k=1}^{N_s}$. Aiming at approximating the solution of the binary infinite horizon optimal control problem \eqref{ocp}, we consider models within the family of feed-forward neural networks (FNNs), for which the choice of $\bu(\cdot)\in\R^{2d}$ as learning target variable can be suboptimal in terms of goodness of fit of the model. A variety of alternatives has been proposed and compared in literature \citep{WANG98}, \citep{kang_algorithms_2020}, \citep{darbon_overcoming_2020}, and a widely popular choice can be to target the associated scalar \emph{value function} $V_\theta(\cdot)\approx V(\cdot)\in\R$, and then recover the feedback as a function of the gradient of $V_\theta(\cdot)$:
\begin{equation}
    \bu_V(\bx) =  -\frac{R^{-1}B^T\nabla V_{\theta}(\bx)}{2}\,,
\end{equation}
where $\nabla V_\theta(\cdot)$ can be efficiently retrieved via automatic differentiation. 

\subsection{Gradient-augmented supervised learning}
Since the learning final goal is to approximate the feedback law $\bu(\cdot)$, the accuracy of the gradient approximation $\nabla V_\theta(\cdot)$ is fundamental. In this direction,  our training is strengthened thanks to a \emph{gradient-augmented} loss function, accounting for not only the approximation error in the target variable, but also for the discrepancy in terms of its gradient. This requires an enriched data-set for the training phase of the supervised learning procedure, including both the value function associated to the binary infinite horizon problem and its gradient $\mathcal{T} = \{\bx^{(k)},V(\bx^{(k)}),\nabla V(\bx^{(k)})\}_{k=1}^{N_s}$. In particular, at every sampled state $\bx$, once the frozen SDRE has been solved for $\Pi$, we consider the ansatz $V(\bx) = \bx^T\Pi\bx$ and we approximate the space derivative  $\nabla V(\bx) \approx 2\Pi\bx$, obtained by neglecting the state dependency of the SDRE solution. 

Even if this ansatz for $V(\cdot)$ approximates the HJB PDE only around a neighborhood of the origin when applied to nonlinear problems, and the target gradient is not exact, this choice allows to conveniently collect the enriched data-set without any computational cost additional to the solution of the ARE associated to the current state. 


\subsection{Network Architecture}
Feed-forward neural networks approximate a function via a chain of composition of layers $l_1,...,l_M$, consisting of an \emph{activation function} $\sigma(\cdot)$ applied component-wise to a linear combination of the layer input variable: 
\begin{equation}
    l_{m}(\by) = \sigma_m(\bA_m \by + \bb_m)\,.
\end{equation} 
The weight matrices $\{\bA_m,\bb_m\}_{m}=\theta$ are parameters to be optimized during the training phase, so that the associated ANN minimizes a suitable \emph{loss function}. In the gradient-augmented settings, we consider a compromise between a fitting functional and a gradient regulation:
\begin{equation}
    \mathcal{L}oss(V,V_\theta) = \mathcal{L}_2(V,V_\theta) + \mu_{dV}\mathcal{L}_2(\nabla V,\nabla V_\theta)\,,
\end{equation}
where $\mathcal{L}_2(f,f_\theta)$ denotes the \emph{mean squared error} (MSE):
\begin{equation}\label{mse}
    \cL_{2}(f, f_{\theta}) := \frac{1}{N_s}\sum_{k=1}^{N_s}\|f(\bx^{(k)}) - f_{\theta}(\bx^{(k)})\|^2\,.
\end{equation}
The number $M$ of layers, their width (i.e. the number of neurons per layer), the activation functions $\sigma_m(\cdot)$, and the loss weight $\mu_{dV}$ are \emph{hyper-parameters} of the model to be optimally tuned so that the trained model not only reaches a good approximation for the training set $\mathcal{X}_t$, but also generalizes outside the training data.

\section{Controlling opinion dynamics}
The aforementioned methodology has been assessed with a numerical test from \citep{albi_mean_2017} dealing with a high-dimensional consensus problem for first-order opinion dynamics governed by the Sznajd model \citep{Sznajd}. Here, the evolution of the state variables is described through the asymmetric interaction kernel $P(\cdot,\cdot)$ defined as follows:
\begin{equation}\label{P}
    P(x_i,x_j) = \beta (1-x_i^2)\,,\qquad \beta\in\R.
\end{equation}
We limit our state space to samples in $\Omega = [-1,1]$ describing the opinions of a large population of voters between two extremal positions $\{-1,1\}$. The interaction kernel models the propensity of agents to change their opinions when interacting with others: the more the agent's opinion is close to the boundary of the domain $\Omega$, the less they are going to influence their peers. A choice of a parameter $\beta < 0$ leads to separation of opinions, meaning that without any external action, the population's opinion is going to concentrate around $x=1$ and $x=-1$ \citep{opinion_models}. {Here, we fix  $\beta = -1$}.

As previously discussed, for a sufficiently high number of agents (here we consider $N_a = 10^5$), the individual-based model \eqref{CStest_dynamics} can be cast in its mean field formulation \eqref{mf}, which is consistent with a kinetic-like equation \eqref{boltzmann} for the evolution in time of the population probability density of having an agent with opinion $x\in\Omega$. This Boltzmann description allows us to approximate the evolution of the mean field dynamics as a limit of binary interactions of sampled couples of agents within the population. 

In order to cast the problem under consideration in semilinear form, we consider the following change of variables 
\begin{equation}\label{cov}
    (x_i,x_j)^{\top}\mapsto(x_i,\Bar{x})
\end{equation}
This allows us to write 
the cost functional \eqref{CStest_cost} in quadratic form \eqref{ocp} with weights $\bR = \gamma/2$ and $\bQ = 2\mathbb{I}_{2}-\mathbb{J}_{2}$, where respectively $\mathbb{I}_{2}$ is the identity, and $\mathbb{J}_{2}$ denotes the matrix full of ones in $\R^{2\times2}$. Similarly, the binary dynamics \eqref{CStest_dynamics} with Sznajd kernel \eqref{P}
\begin{equation}\dot{\bx}=
    \begin{cases}
    \dot{x}_i = \dfrac{\beta}{2}(1-x_i^2)(x_j-x_i) + u_i\\[0.5em]
    \dot{x}_j = \dfrac{\beta}{2}(1-x_j^2)(x_i-x_j) + u_j
    \end{cases}
\end{equation}
can be written, after the change of variable \eqref{cov}, in semilinear form as
\begin{equation}
   \begin{bmatrix}
    x_i\\
    \Bar{x}
    \end{bmatrix} = 
    \begin{bmatrix}
    -P(x_i,\Bar{x}) & P(x_i,\Bar{x})\\
    -\Bar{P}(x_i,\Bar{x}) & \Bar{P}(x_i,\Bar{x})
    \end{bmatrix}\begin{bmatrix}
    x_i\\
    \Bar{x}
    \end{bmatrix} + \mathbb{I}_{2}\bu
\end{equation}
where
\begin{equation*}
    P(x_i,\Bar{x}) = \beta (1-x_i^2)\qquad \Bar{P}(x_i,\Bar{x}) = \beta\big((2\Bar{x}-x_i)^2-x_i^2\big).
\end{equation*}

For populating the dataset $\mathcal{X}_t$ we uniformly sample $N_s = 10^3$ states $\bx^{(k)}=(x_i,\Bar{x})^{(k)}$ and we apply \eqref{cov}. For every sampled current state of the system, we compute the state-dependent SDRE coefficients and we rely on the \texttt{lqr} routine in MATLAB for solving the associated ARE for $\Pi(\bx^{(k)})$. With this suboptimal SDRE solution, the training set $\mathcal{T}_V=\{\bx^{(k)},V^{(k)},\nabla V^{(k)}\}_{k=1}^{Ns}$ is computed with $V^{(k)} = \bx^{(k)T}\Pi(\bx^{(k)})\bx^{(k)}$, and $\nabla V^{(k)} = 2\Pi(\bx^{(k)})\bx^{(k)}$. An additional dataset is generated, with the purpose of comparing the gradient-augmented approximation $\bu_V$ with the direct approximation of the feedback law $\bu_\theta$: $\mathcal{T}_u=\{\bx^{(k)},u^{(k)}\}_{k=1}^{Ns}$, where $\bu^{(k)} = -\bR^{-1} \bB^{T} \Pi(\bx^{(k)})\bx^{(k)}$. For the training of the model $\bu_\theta$ the loss function was the MSE \eqref{mse}.

Once the datasets $\mathcal{T}_V,\,\mathcal{T}_u$ have been computed, they have been split into \emph{training sets} and \emph{validation sets}, with a ratio of $80/20$. {The ANN architectures for both the ANN $\bu_\theta$ and $V_\theta$  (together with the loss weight $\mu_{dV}\in[0,2]$) were chosen accordingly with the goodness of fit of the model evaluated within the validation samples: this hyper-parameter tuning phase has been dealt with via a grid search in the parameter space by maximizing the precision of the trained model, by means of minimization of the \emph{mean relative error} (MRE).}

The desired architecture was identified in both cases to be a FFN with $M=4$, having identity activation function for the input and output layers $\sigma_{1,4}(\by)=\by$, and a sigmoid function for the remaining ones $\sigma_{2,3}(\by)=(1+e^{-\by})^{-1}$. The hidden layers were populated by $100$ neurons per each, while the dimension of the state space $\Omega$, and the scalar nature of the target of $V_\theta$ defined the depths of the remaining layers. For $V_\theta$, the best configuration of hyper-parameters set the loss weight to $\mu_{dV}=0.05$.

The goodness of fit of the trained models is finally evaluated in a \emph{test set}, a uniform grid of $N_v=10^5$ points within the state space, where the approximated control is compared with the pointwise computation through the SDRE solution. Goodness of fit of trained models in both tests are presented in Table \ref{tab:error_test}.

\begin{table}[h]
    \centering
    \caption{Goodness of fit in terms of: MSE, \emph{coefficient of determination} $r^2$, and MRE.}
    \begin{tabular}{c l l l l}
         \emph{target} &  & $MSE$ & $r^2$ &  $MRE$ \\[2pt]
        \hline
        
        \multirow{2}{*}{$V_\theta$} 
        & \multicolumn{1}{l}{$\mu_{dv} = 0.05$} & \multicolumn{1}{l}{$9.24\times10^{-6}$} & \multicolumn{1}{l}{$0.96480$} & \multicolumn{1}{l}{$0.0195$}\\
        & \multicolumn{1}{l}{$\mu_{dv} = 0$} & \multicolumn{1}{l}{$3.50\times10^{-5}$} & \multicolumn{1}{l}{$0.86674$} & \multicolumn{1}{l}{$0.2779$}\\[5pt]
        
        \multirow{2}{*}{$dV_\theta$} 
        & \multicolumn{1}{l}{$\mu_{dv} = 0.05$} & \multicolumn{1}{l}{$3.71\times10^{-7}$} & \multicolumn{1}{l}{$0.99992$}& \multicolumn{1}{l}{$0.0079$}\\
        & \multicolumn{1}{l}{$\mu_{dv} = 0$} & \multicolumn{1}{l}{$4.48\times10^{-5}$} & \multicolumn{1}{l}{$0.99017$} & \multicolumn{1}{l}{$0.0987$}\\[5pt]
        
        \multirow{2}{*}{$u_V$} 
        & \multicolumn{1}{l}{$\mu_{dv} = 0.05$} & \multicolumn{1}{l}{$5.94\times10^{-4}$} & \multicolumn{1}{l}{$0.99992$} & \multicolumn{1}{l}{$0.0079$}\\
        & \multicolumn{1}{l}{$\mu_{dv} = 0$} & \multicolumn{1}{l}{$0.071668$} & \multicolumn{1}{l}{$0.98417$} & \multicolumn{1}{l}{$0.0987$}\\[5pt]
        $u_\theta$ & & $0.002778$ & $0.99962$ & $0.0195$\\
    \end{tabular}
    \centering
    \label{tab:error_test}
\end{table}

\begin{figure}[ht]
	\centering
    \includegraphics[width=0.46 \textwidth]{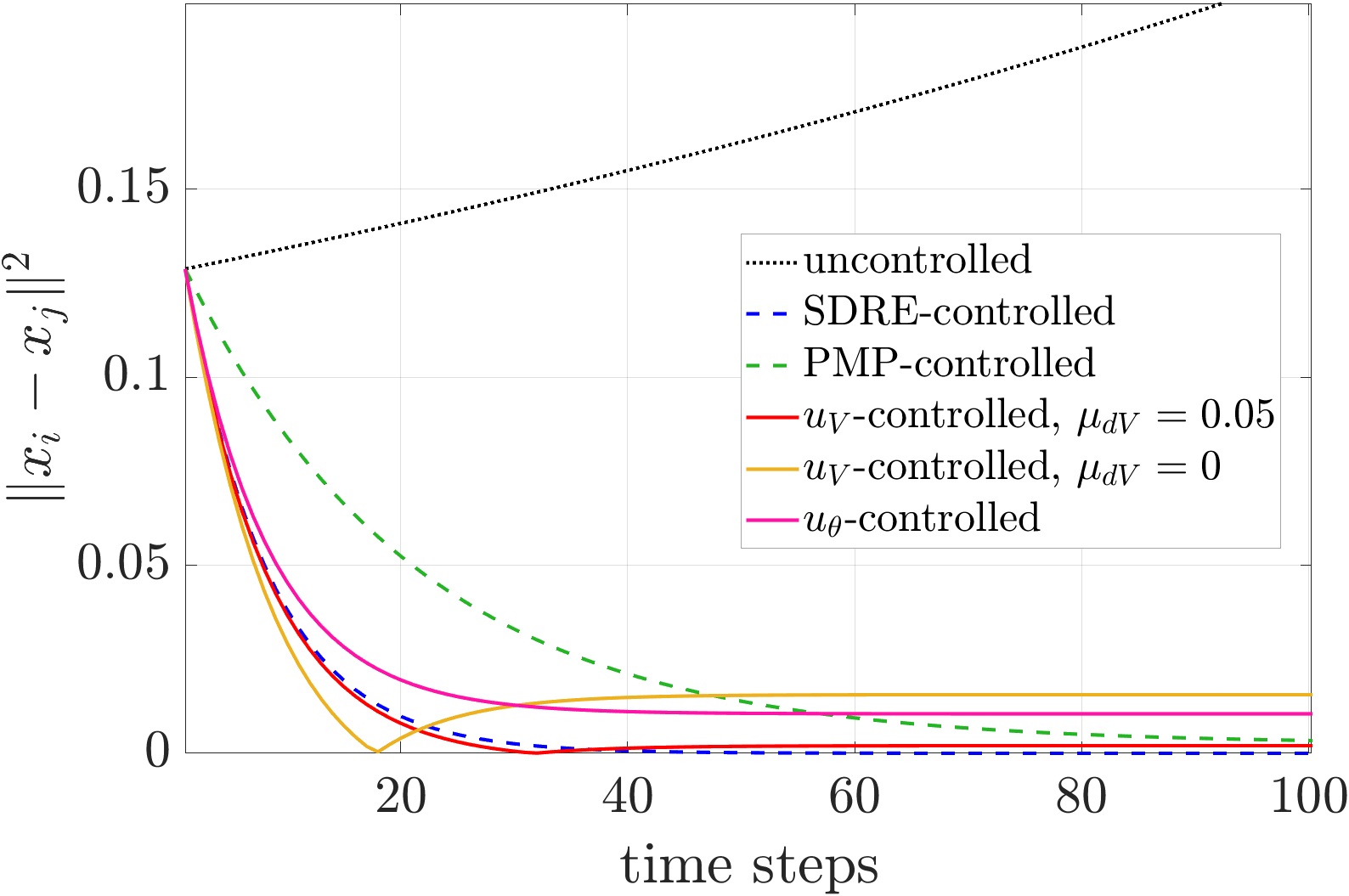}
	\caption{Evolution of the Euclidean distance within a sampled couple of agents in both uncontrolled ($u(x)\equiv0\,\forall x\in\Omega$) and controlled settings. {The SDRE feedback law $\bu$ succeeds at steering the couple system towards consensus much faster than the open-loop control variable (for which we need $T>100$ to reach consensus).} Among the different approximations, the feedback $\bu_V$, resulting from the gradient-augmented model, leads to the best performance.}
	\label{t2traj}
\end{figure}

With the trained models $V_\theta$ and $\bu_\theta$ for the binary interactions, we compare the evolution of a sampled couple under the action of the different controls: the suboptimal feedback law $\bu$ obtained with the SDRE approach and its approximations $\bu_V$ and $\bu_\theta$. A further comparison is done w.r.t. the open-loop control variable obtained by solving Pontryagin's optimality conditions (PMP) holding in finite horizon settings. Aiming at approaching the feedback behaviour in PMP settings, we consider a time horizon $T$ large enough for the system to reach consensus. A plot of the dynamics of a sampled couple of agents  $x_i,x_j\in\Omega = [-1,1]$ can be seen in Figure \ref{t2traj}. 

\begin{figure}[ht]
	\centering
    \includegraphics[width=0.4 \textwidth]{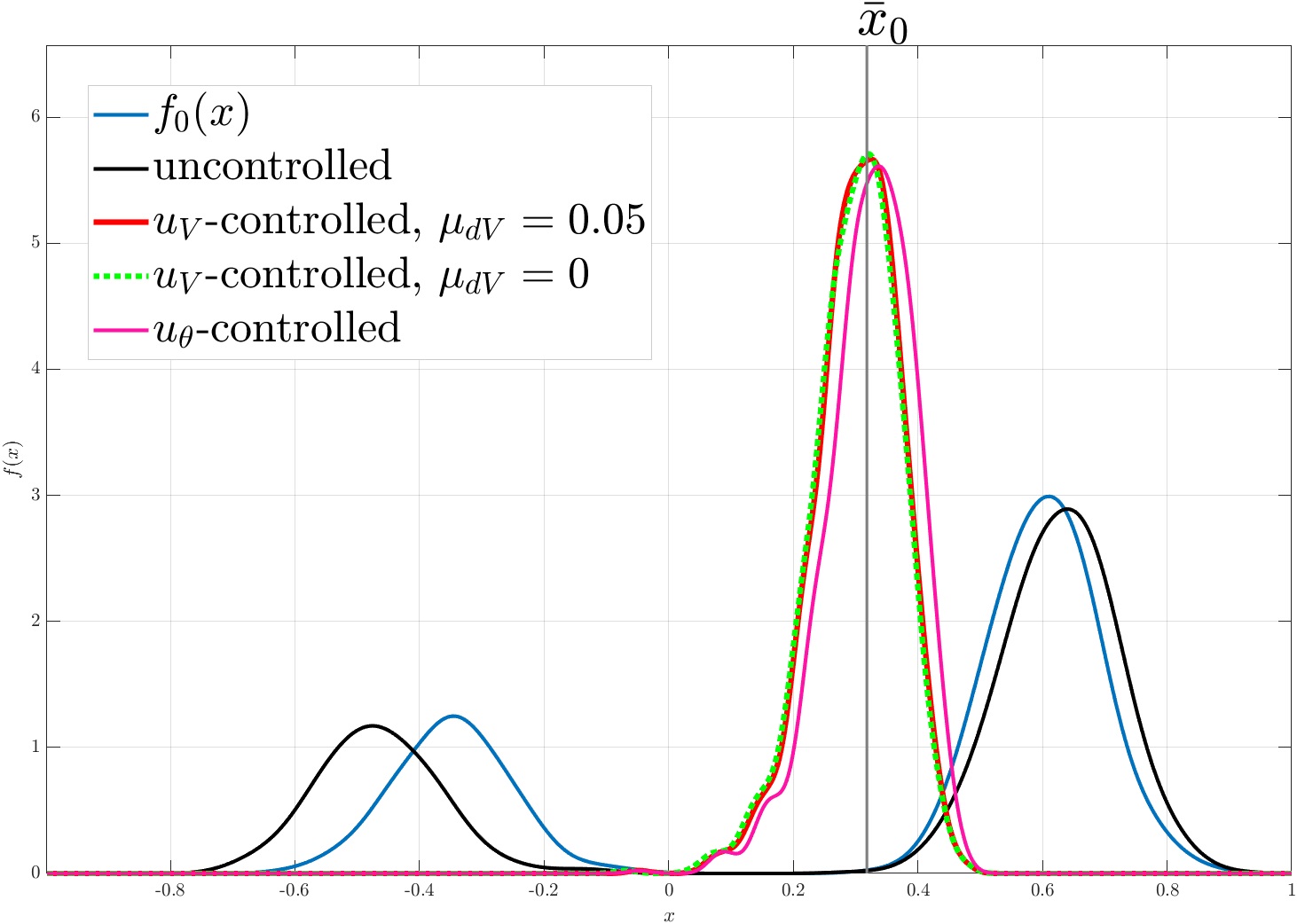}
	\caption{Comparison between the initial density function $f_0(x)$ with the sampled probability density obtained as limit of binary interactions controlled by approximated control variables $\bu_V$ and $\bu_\theta$. The controlled system has only been evolved for 10 time steps and yet it is already concentrating around the target opinion (consensus) for both the standard and the gradient-augmented $\bu_V$ feedback laws. The action of $\bu_\theta$ can be seen to steer the agents towards a slightly different configuration. The uncontrolled dynamics are leading to opinion separation, consistently with the parameter choice $\beta = -1<0$.}
	\label{traj_compare}
\end{figure}

\begin{figure}[h]
	\centering
    \includegraphics[width=0.4 \textwidth]{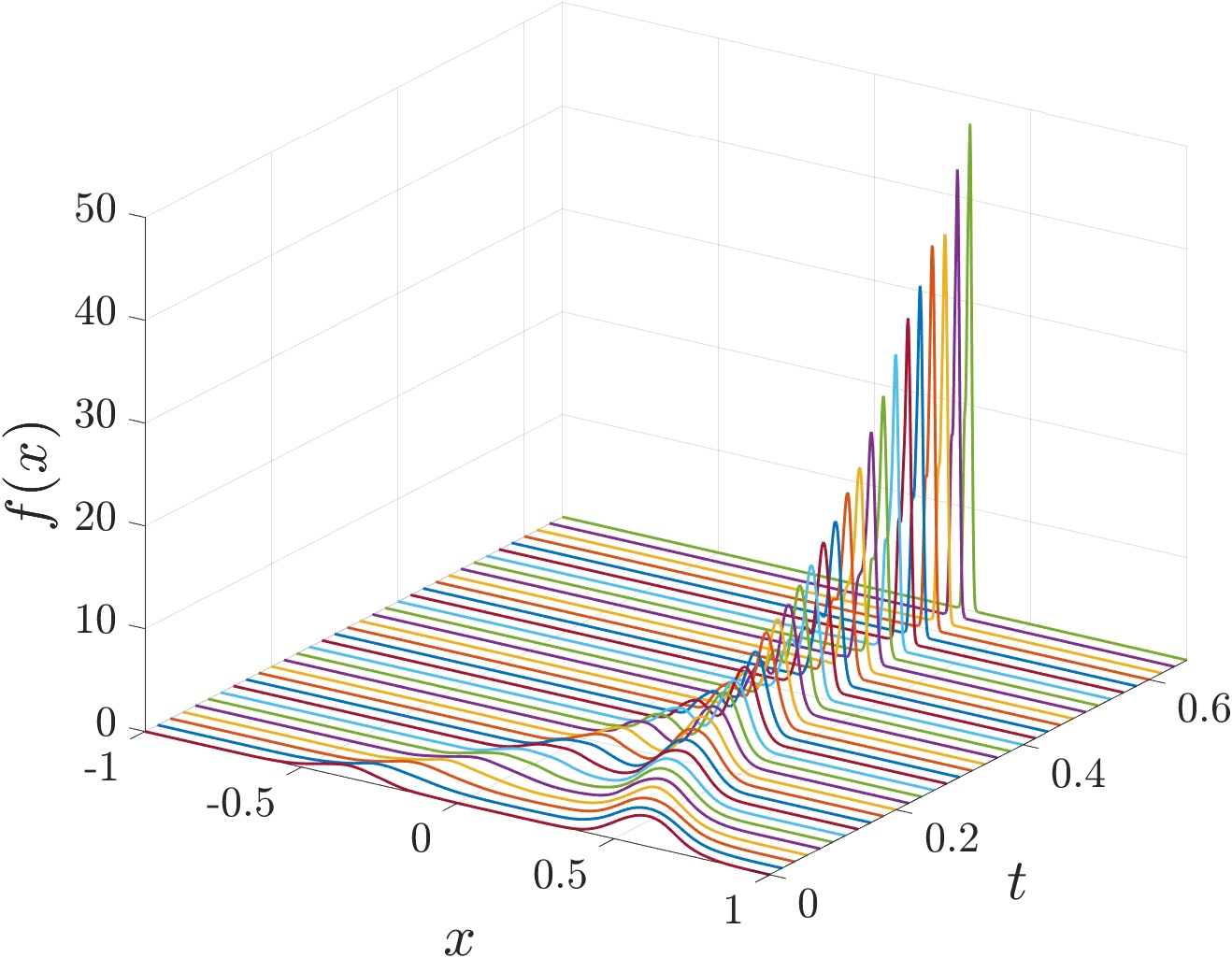}
	\caption{Evolution of the $\bu_V$-controlled probability density of finding an agent with opinion $x$ versus time (seconds). From a double-peaked density obtained as a mixture of normal densities, the opinions of the agents rapidly converge towards consensus.}
	\label{traj_uV}
\end{figure}

Finally, we can rely on the binary interactions controlled via $\bu_V,\,\bu_\theta$ in order to approximate the behaviour of the whole population. In particular, with the choice of time-step $\Delta t = \varepsilon$ in \eqref{euler}, we allow each one of the agents to interact with someone else at every update. This means that at each time step, we can sample $N_a/2$ couples within the population and then act on their interactions by means of a feedback variable. Every couple evolves according to a forward Euler scheme with time-step $\Delta t$, after which the population density function is updated to be the sampled density of all the agents. In this way we approximate the behaviour of the controlled population from a mean field viewpoint, by only solving many $2$-dimensional infinite horizon optimal control problems. The time evolution of the population probability density function influenced through $\bu_V$ can be seen in Figure  \ref{traj_uV}. In Figure \ref{traj_compare}, we compare the given initial distribution $f_0(x)$ with its time evolution according to controlled mean field dynamics by means of a variety of feedback laws.

\section{Conclusions}
In this paper a mixture of approximation techniques for solving optimal control of multi-agent systems has been discussed and numerically tested. The first approximation step coincides with considering a mean field formulation of the controlled dynamics, so that the number of agents populating the system no longer contributes to the dimensionality of the problem. Then, the complexity of the solution of such a mean field optimal control problem has been further reduced thanks to a description of the population dynamics from a kinetic viewpoint, by means of a Bolzmann equation for the time evolution of the population density. This formulation has the advantage that the complexity of the associated solution is dramatically reduced with respect to the mean field optimal control, still retaining the ability to influence the population as a whole. Finally, a gradient-augmented supervised learning model has been trained for approximating the suboptimal SDRE solution of the reduced Bolzmann binary interactions. A comparison between the proposed model and the direct approximation of the feedback law in a numerical example has highlighted an outstanding performance of the former. In the future we will assess the proposed methodology in higher dimensional dynamics, where the supervised learning of the feedback map is essential to enable computational feasibility of the kinetic approach. 
\bibliography{bibliography}
\end{document}